\documentclass[11pt,twoside]{article}
\usepackage{times}
\usepackage{amsmath,amssymb}
\usepackage{color}
\usepackage{graphicx}

\allowdisplaybreaks

\newcommand{\pf}{\noindent {\bf Proof. \hspace{2mm}}}
\newcommand{\ef}{ \hfill $ \blacksquare $ \vskip 3mm}
\newcommand{\be}{\begin{equation}}
\newcommand{\ee}{\end{equation}}
\newcommand{\bea}{\begin{eqnarray}}
\newcommand{\eea}{\end{eqnarray}}

\def\nn{\nonumber}

\def\p{\partial}
\def\ve{\varepsilon}
\def\f{\frac}

\def\O{\Omega}
\def\th{\theta}
\def\g{\gamma}

\def\p{\partial}

\def\ve{\varepsilon}
\def\f{\frac}

\def\th{\theta}

\def\g{\gamma}

\def\ds{\displaystyle}

\def\th{\theta}

\def\g{\gamma}

\def\nn{\nonumber}

\def\p{\partial}
\def\ve{\varepsilon}
\def\f{\frac}

\def\O{\Omega}
\def\th{\theta}
\def\g{\gamma}

\def\p{\partial}

\def\ve{\varepsilon}
\def\f{\frac}

\def\th{\theta}

\def\g{\gamma}

\def\ds{\displaystyle}

\def\th{\theta}

\def\g{\gamma}

 \allowdisplaybreaks

\begin{document}
 \footskip=0pt
 \footnotesep=2pt
\let\oldsection\section
\renewcommand\section{\setcounter{equation}{0}\oldsection}
\renewcommand\thesection{\arabic{section}}
\renewcommand\theequation{\thesection.\arabic{equation}}
\newtheorem{claim}{\noindent Claim}[section]
\newtheorem{theorem}{\noindent Theorem}[section]
\newtheorem{lemma}{\noindent Lemma}[section]
\newtheorem{proposition}{\noindent Proposition}[section]
\newtheorem{definition}{\noindent Definition}[section]
\newtheorem{remark}{\noindent Remark}[section]
\newtheorem{corollary}{\noindent Corollary}[section]
\newtheorem{example}{\noindent Example}[section]

\title{On the global existence and stability of 3-D viscous cylindrical circulatory  flows}
\author{ }

\author{Yin,
Huicheng$^{1,*}$; \quad Zhang, Lin$^{2}$\footnote{*Yin Huicheng (huicheng$@$nju.edu.cn, 05407@njnu.edu.cn) and Zhang Lin (lynzhung$@$gmail.com) were
supported by the NSFC (No.~11025105) and A Project Funded by the Priority Academic Program Development of
Jiangsu Higher Education Institutions.}\vspace{0.5cm}\\
\small 1. School of Mathematical Sciences, Jiangsu Provincial Key Laboratory for Numerical Simulation\\
\small of Large Scale Complex Systems, Nanjing Normal University, Nanjing 210023, China.\\
\small 2. Department of Mathematics and IMS, Nanjing University, Nanjing 210093, China.\\
}

\date{}
\maketitle

\vskip 0.3 true cm

\centerline {\bf Abstract} \vskip 0.3 true cm

In this paper, we are concerned with the
global existence and stability of a 3-D perturbed viscous circulatory
flow around an infinite long cylinder. This flow is described by
3-D compressible Navier-Stokes equations.
By introducing some suitably weighted energy spaces and
establishing a priori estimates, we show that the
3-D cylindrical symmetric circulatory flow is globally stable in time when
the corresponding initial states are perturbed suitably small.

\vskip 0.3 true cm

 \vskip 0.3 true cm


\section{Introduction}

In this paper, we are concerned with the global stability problem of cylindrical symmetric circulatory
flows for the three-dimensional compressible Navier-Stokes equations (see Figure 1 below). The
compressible Navier-Stokes equations in three space dimensions are
\begin{equation}
\left\{
\begin{aligned}
&\p_t\rho +div(\rho u)=0, \\
&\rho \p_tu +\rho u\cdot \nabla u +\nabla P(\rho) = \nu_1 \Delta u+\nu_2\nabla div u,
\end{aligned}
\right.
\end{equation}
where $\rho>0$ is the density, $u=(u_1, u_2, u_3)$ is the velocity,
$\nu_1> 0$ and $\nu_1+\nu_2>0$ hold, and the state equation is given by $P(\rho)=A\rho^\gamma$ with
the constants $A>0$ and $\gamma>1$.

We now give a mathematical description on the 3-D viscous cylindrical
flow around an infinite long cylinder $\{x=(x_1,x_2, z)\in\Bbb R^3: r=\sqrt{x_1^2+x_2^2}\le 1, z\in\Bbb R\}$.
Set $\O=\{(r,z): r>1, z\in\Bbb R\}$ and $(\rho(t,x), u(t,x))=(\rho(t,r, z),
u_r(t, r, z)\ds\f{x'}{r}+u_{\th}(t, r, z)\f{{x'}^{\bot}}{r}, u_z(t,r,z))$, where $x'=(x_1, x_2)$ and ${x'}^{\bot}=(-x_2, x_1)$.
In this case, (1.1) has the following equivalent form in $[0, \infty)\times \O$:
\bea
&&\p_t \rho +\f{1}{r}\p_r(r\rho u_r)+\p_z(\rho u_z)=0,  \label{cc1.1} \\
&&\rho \p_t u_r +\rho (u_r \p_r u_r+u_z\p_zu_r -\f{u_\theta^2}{r})+\p_r P(\rho) \nn\\
&&\quad \quad =\nu_1\biggl(
\p_r(\f{1}{r} \p_r (r u_r))+\p_z^2 u_r\biggl) +\nu_2 \p_r \bigg (\f{1}{r}\p_r(r  u_r)+\p_zu_z\bigg),  \label{cc1.2}\\
&& \rho \p_t u_\theta +\rho ( u_r \p_r u_\theta+u_z\p_z u_\theta +\f{u_\theta u_r}{r})=\nu_1 \biggl(
\p_r(\f{1}{r} \p_r (r u_\theta))+\p_z^2 u_\theta\biggl), \label{cc1.3}\\
&&
\rho \p_t u_z +\rho(u_r \p_r u_z+u_z\p_zu_z)+\p_z P(\rho)\nn \\
&&\quad \quad =\nu_1\biggl(
\p_r^2u_z+\p_z^2 u_z+\f{1}{r}\p_r u_z\biggl) +\nu_2\p_z\bigg(\f{1}{r}\p_r(r  u_r)+\p_zu_z\bigg). \label{cc1.4}
\eea

\begin{figure}[htbp]
\centering\includegraphics[width=10cm,height=6.5cm]{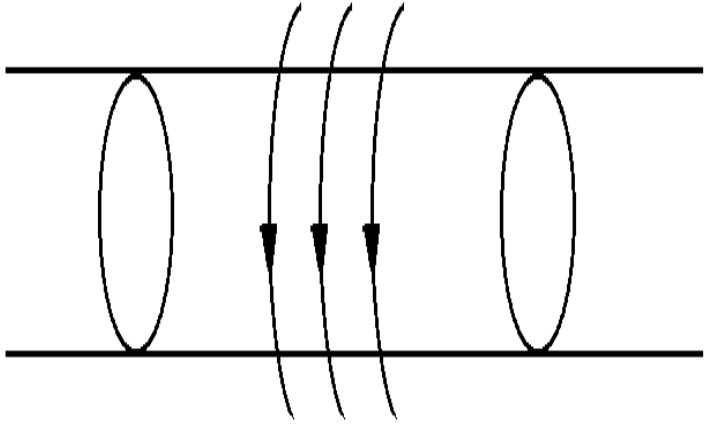}
\caption{Cylindrical circulatory flow around a cylinder}\label{fig,1}
\end{figure}

\vskip 0.4 true cm

We start to look for a special steady circulatory flow $(\bar\rho(r), \bar u_r(r), \bar u_\theta(r), 0)$
of (1.2)-(1.5) around the cylinder $\{(r,z): r\le 1, z\in\Bbb R\}$. Such a flow is called a background solution of (1.2)-(1.5)
in the whole paper. By (1.2)-(1.5), one knows that $(\bar\rho(r), \bar u_r(r), \bar u_\theta(r))$
satisfies
\begin{equation}
\left\{
\begin{aligned}
&\f{1}{r}\p_r(r\bar\rho \bar u_r)=0, \\
&\bar \rho (\bar u_r \p_r \bar u_r -\f{\bar u_\theta^2}{r})+\p_r P(\bar \rho) =(\nu_1+\nu_2)
\p_r(\f{1}{r} \p_r (r \bar u_r)), \\
&\bar \rho (\bar  u_r \p_r \bar u_\theta +\f{\bar u_\theta \bar u_r}{r})=\nu_1
\p_r(\f{1}{r} \p_r (r \bar u_\theta)).
\end{aligned}
\right.\label{O}
\end{equation}
On the other hand, in order to solve (1.6), one naturally poses a boundary condition
on $\Sigma'=\{x': r=1\}$ as follows
\bea
&&(\bar\rho, \bar u_r, \bar u_\th)|_{r=1}=(\bar\rho_0, 0, M_0),\label{A}
\eea
where $\bar\rho_0>0$ and $M_0>0$ are given constants.
In addition, one also requires that at infinity
\bea
&&\ds\lim_{r\to\infty}\bar u_r(r)=0,\quad \ds\lim_{r\to\infty}\bar  u_\th(r)=0
,\quad \ds\lim_{r\to\infty}\bar  u_z(r)=0. \label{B}
\eea
As illustrated in [3] or [22], it is easy to know that (1.6) with (1.7)-(1.8) has a unique solution for $r\ge 1$
\bea
&&(\bar\rho(r), \bar u_r(r), \bar u_\theta(r))=\biggl(\bigl(\bar\rho_0^{\g-1}
+\ds\f{(\g-1)M_0^2}{2A\g}(1-\f{1}{r^2})\bigr)^{\f{1}{\g-1}},\quad 0,\quad \ds\f{M_0}{r}\biggr).
\label{C}
\eea

In the paper, we focus on the global stability problem of the background solution $(\bar\rho(r), \bar u_r(r), \bar u_\theta(r),$
$0)$.
Namely, the global
solution problem of (1.2)-(1.5) in the domain $[0, \infty)\times\O$ will be studied under the
following perturbed initial value conditions:
\bea
&&(\rho, u_r, u_\theta, u_z)(0,r,z)=(\bar\rho(r)+\rho_0(r,z),u_r^0(r,z),
\bar u_\theta(r)+u_\theta^0(r,z),u_z^0(r,z)),\\
&&\ds\lim_{r^2+z^2\to\infty}u_r(t, r,z)=\ds\lim_{r^2+z^2\to\infty}u_\th(t, r,z)=\ds\lim_{r^2+z^2\to\infty}u_z(t, r,z)=0,\\
&&(u_r(t, r,z), u_\th(t, r,z), u_z(t, r,z))|_{r=1}=(0, M_0, 0),
\label{D}
\eea
where $\rho_0\in H_0^2(\O)$ and $(u_r^0, u_\theta^0, u_z^0)\in H^3_0(\O)$.

Let
$\ds\rho(t,r,z) =\bar\rho(r) +\phi(t,r,z)$, $\ds u_r(t,r,z)=v_r(t,r,z)$, $\ds u_\theta(t,r,z) =\bar u_{\th}(r)+v_\theta(t,r,z)$, $\ds u_z(t,r,z)=v_z(t,r,z)$, and $v=(v_r, v_\theta, v_z)$. Then equations \eqref{cc1.1}-\eqref{cc1.4} together with  (1.10)-(1.12) can be written as
\bea
&&\p_t \phi +\f{1}{r}\p_r (r \bar\rho v_r)+\p_z(\bar\rho v_z)=f, \label{ccb1.10}\\
&&\p_t v_r -\f{2M_0}{r^2}v_\theta +\gamma \p_r({\bar\rho}^{\gamma-2}  \phi)
-\f{\nu_1}{\bar\rho}\bigg(\p_r(\f{1}{r}\p_r(rv_r))+\p_z^2 v_r\biggl)-\f{\nu_2}{\bar\rho}\p_r\bigg(\f{1}{r}\p_r (r  v_r)
+\p_z  v_z\bigg)=g_1, \nn\\
&&\label{cc1.11}\\
&&\p_t v_\theta
-\f{\nu_1}{\bar\rho}\biggl(\p_r(\f{1}{r}\p_r(rv_\theta))+\p_z^2 v_\theta\biggl)=g_2,\label{cc1.7} \\
&&\p_t v_z+\gamma \p_z ({\bar\rho}^{\gamma-2} \phi)-\f{\nu_1}{\bar\rho}\biggl( \p_r^2 v_z +\p_z^2 v_z +\f{1}{r}\p_r v_z \biggl)-\f{\nu_2}{\bar\rho}\p_z\bigg(\f{1}{r}\p_r (r  v_r)+\p_z v_z\bigg)=g_3\label{cc1.8}
\eea
with the initial-boundary value conditions
\bea
&&(\rho, v)|_{t=0}=(\phi_0,v_0)\equiv(\rho_0(r,z),u_r^0(r,z),
u_\theta^0(r,z),u_z^0(r,z)), \label{cc1.10}\\
&&v|_{r=1}=(0,0,0), \qquad \ds\lim_{r^2+z^2 \to +\infty}v=(0,0,0),\label{cc1.9}
\eea
where
\bea
&& f= -\f{1}{r}\p_r (r\phi v_r)-\p_z(\phi v_z), \nn \\
&&g_1=\f{v_\theta^2}{r}-v_r\p_r v_r-v_z\p_z v_r -\p_r Q(\bar\rho,\phi)
-\f{\nu_1 \phi}{(\phi+\bar\rho)\bar\rho}\biggl(\p_r(\f{1}{r}\p_r(rv_r))+\p_z^2 v_r \biggl)\nn\\
&&\quad \quad -\f{\nu_2 \phi}{(\phi+\bar\rho)\bar\rho}\p_r\bigg(\f{1}{r}\p_r (r  v_r)+\p_z v_z\bigg),\nn\\
&&g_2 =-v_r \p_r v_\theta -v_z\p_z v_\theta -\f{v_\theta v_r}{r}
-\f{\nu_1 \phi}{(\phi+\bar\rho)\bar\rho}\biggl(\p_r(\f{1}{r}\p_r(rv_\theta))+\p_z^2 v_\theta\biggl), \nn\\
&&g_3=-v_r\p_r v_z-v_z\p_z v_z  -\p_z Q(\bar\rho,\phi) -\f{\nu_1 \phi}{(\phi+\bar\rho)\bar\rho}\biggl(\p_r^2 v_z+\p_z^2 v_z +\f{1}{r} \p_r v_r \biggl)\nn\\
&&\quad \quad -\f{\nu_2 \phi}{(\phi+\bar\rho)\bar\rho}\p_z\bigg(\f{1}{r}\p_r (r  v_r)+\p_z v_z\bigg),\nn
\eea
and
$$
Q(\bar\rho,\phi)=\f{\gamma(\gamma-2)}{2} \phi^2 \int_0^1 (\bar\rho+s\phi)^{\gamma-3} sds.
$$

To state our  main results conveniently, we now introduce the following notations:
for $w_1, w_2\in L^2(\Omega)$, set
$$
(w_1, w_2)=\int_\Omega w_1w_2 drdz.
$$
In addition, $w\in L^p(\Omega) $ $(1\le p < \infty)$  means that
$$\|w\|_{L^p}=\|w\|_p=\bigg(\int_\Omega |w|^p  drdz\bigg)^\f{1}{p}< +\infty.$$
And define
$$
L_r^p(\Omega)=\{w\in \mathcal{D}'(\Omega): \|w\|_{L_r^p}< +\infty \},
$$
where
$$
\|w\|_{L_r^p}=\bigg(\int_\Omega |w|^p rdrdz\bigg)^\f{1}{p}.
$$
Set $D=(\p_r,\p_z)$ and  define for $k\in\Bbb N\cup\{0\}$
$$
\|w\|_{\tilde H^k}^2=\sum \limits_{j=0}^k  \|\sqrt{r}D^j w\|_2^2.
$$
Denote by
$$
\tilde H^k =\{w\in \mathcal{D}'(\Omega): \|w\|_{\tilde H^k}<+\infty\}.
$$
The main conclusion in the paper is:
\begin{theorem}
There exists a constant $\varepsilon>0$ such that if  $\|\phi_0\|_{\tilde H^2}+\|v_0\|_{\tilde H^3} \le \varepsilon$, then
problem (1.13)-(1.16) together with (1.17)-(1.18)  has a unique global solution $(\phi, v)\in C([0,\infty),\tilde H^2\times\tilde H^3)$
satisfying
$$
\|\phi\|_{\tilde{H^2}}^2+\|v\|_{\tilde{H^3}}^2+\int_0^{\infty}(\|D(\bar\rho^{\gamma-2}\phi)\|_{\tilde H^1}^2+
\|D v\|_{\tilde H^3}^2)d\tau \le C (\|\phi_0\|_{\tilde H^2}^2+\|v_0\|_{\tilde H^3}^2).
$$
\end{theorem}

{\bf Remark 1.1.} {\it For the original problem (1.2)-(1.5) with (1.10)-(1.12), one knows from Theorem 1.1 that the
perturbed cylindrical symmetric circulatory flows are globally stable.}

\vskip 0.1 true cm
{\bf Remark 1.2.} {\it So far there have been extensive results on the global spherically symmetric
(or helically symmetric) weak/strong /classical solutions to the
compressible Navier-Stokes equations (in this case, the solution admits a form $(\rho(t,x), u(t,x))$ $=(\rho(t,r), U(t,r)\ds\f{x}{r})$
or $(\rho(t,x), u(t,x))=(\rho(t,r, z),
u_r(t, r, z)\ds\f{x'}{r}+u_{\th}(t, r, z)\f{{x'}^{\bot}}{r}, u_z(t,r,z))$ with $x'=(x_1, x_2)$ and ${x'}^{\bot}=(-x_2, x_1)$),
one can see [2], [4], [7], [14-15], [19], [21] and the references therein. Here, we point out that our initial data
in (1.10) has no bounded energy (due to $u(0,x)\sim\ds\f{1}{r}$), which is a little different from the cases in the aforementioned
references.}

\vskip 0.1 true cm
{\bf Remark 1.3.} {\it  For the Cauchy problem or initial-boundary value problem in exterior domain
of 3-D compressible Navier-Stokes equations, when the initial data are in some suitably weighted
energy spaces or are of  small perturbations with respect to the constant states, many authors
have established the local/global existence of weak/strong/classical solutions in appropriate function spaces,
one can find the details in  [1], [5-6], [8], [13], [16-18], [20] and so on. If we intend to study
the general (not cylindrical symmetric) global perturbation problem of 3-D circulatory flows for (1.1), the methods applied in the
above references cannot be applied directly since our perturbed initial data are different from those
(for examples, our initial data have not finite energies or are not of the small perturbations of constant states).
On the other hand, motivated by the results and methods in [9-10] and [11-12], where
the global stabilities and large time behaviors of the perturbed constant equilibrium on the half
space, and  of the perturbed plane Couette flow are studied respectively
when the Reynolds and Mach numbers are sufficiently small, we hope that the global stability of generally perturbed viscous circulatory flows
can be established in our future research.}

Let's recall some previous works which are related to our results. For the initial-boundary value problem of (1.1),
the local classical solution is obtained in [20] with $\rho_0$ being positive and bounded. Applying the energy methods in Sobolev spaces, the authors
in [17] established  the global existence of classical solutions to (1.1) when the
initial data are of small perturbations for a non-vacuum constant state and no slip boundary conditions
are posed. Recently, for the case that the initial density is allowed to vanish and even has compact support
and the smooth initial data are of small total energy,
the authors in [8] established the  global existence and uniqueness of
classical solutions whose corresponding far fields are vacuum or non-vacuum.
For the arbitrary initial data with finite total energies,
the global existence of weak solutions to 3-D compressible Navier-Stokes equations has been established by P. L. Lions
in [13] for suitably large adiabatic exponent $\g$, and subsequently this result  was improved
to the cases of $\g>\f32$ for the general solutions in [5]
and $\g>1$ for cylindrically symmetric solutions in [15] respectively.
In addition, D. Hoff in [7]
showed the existence of spherically symmetric  weak solutions for $\g= 1$ and
discontinuous initial data. Here we point out that the corresponding background solution of
(1.2)-(1.5) is not a constant state, which is different from those situations in the
aforementioned references; in addition, compared with reference [22], where
the global existence and stability of a 2-D perturbed viscous symmetric circulatory
flow around a disc are established, the analysis  in the present paper is more
involved due to the multi-dimensional spaces.

To prove Theorem 1.1,  we require to establish some global weighted energy estimates of the solution  $(\phi, v)$
to (1.13)-(1.16). Thanks to delicate analysis, the uniform weighted estimates
of $(\phi, v)$ are obtained by making full use of the properties (for instance, $\p_r\bar\rho\sim\ds\f{1}{r^3}$)
of the background solution
and choosing suitable multipliers.
Based on this and the local existence result of classical solution to (1.13)-(1.16) with (1.17)-(1.18),
Theorem 1.1 is shown by the continuity argument.

The paper is organized as follows: In \S2, we derive some uniform  energy estimates from
the linearized parts of (1.13)-(1.16). From this,
some uniform weighted energy inequalities of  $(\phi, v)$
are obtained and subsequently the proof of Theorem 1.1 is completed in $\S 3$.

\section{Some Elementary Estimates}

In this section, we establish some basic weighted energy inequalities on the solution  $(\phi, v)$ of (1.13)
-(1.16) with (1.17)-(1.18).

\begin{lemma}\label{pcl1} ({\bf Weighted $L^2-$estimate of $(\phi, v)$}). For the solution
$(\phi, v)\in C([0,\infty),\tilde H^2\times\tilde H^3)$ of problem (1.13)-(1.16) with (1.17)-(1.18),
we have
\bea
&& \|\sqrt{\bar\rho r}v \|_2^2 +\|\sqrt{{\bar\rho}^{\gamma-2} r}\phi\|_2^2+ \int_0^t\biggl( \|\f{v_r}{\sqrt{r}}\|_2^2+\|\f{v_\theta}{\sqrt{r}}\|_2^2+\|\sqrt{r}D v\|_2^2
+\|\sqrt{r}(\f{1}{r}\p_r(rv_r)+\p_z v_z)\|_2^2\biggl)d\tau\nn\\
&&\le C(\|\sqrt{ r}v_0 \|_2^2 +\|\sqrt{r}\phi_0\|_2^2)+ C \int_0^t A_1d\tau, \nn
\eea
where and below $C>0$ stands for a generic constant, and $A_1=|(g, \bar\rho  r v)|+|(f, {\bar\rho}^{\gamma-2} r\phi)|$
with $g=(g_1, g_2, g_3)$.
\end{lemma}
\pf  It follows from
$\ds\int_\Omega \eqref{cc1.11} \times \bar\rho r v_r drdz$,
$\ds\int_\Omega \eqref{cc1.7} \times \bar\rho r v_\theta drdz$ and  $\ds\int_\Omega \eqref{cc1.8} \times \bar\rho r v_z drdz$  that

\bea
&&\f{1}{2}\f{d}{dt}\|\sqrt{\bar\rho r}v_r\|_2^2
+\gamma (\p_r(\bar\rho^{\gamma-2} \phi), \bar\rho r v_r)
+\nu_1(\|\f{v_r}{\sqrt{r}}\|_2^2+\|\sqrt{r}\p_r v_r\|_2^2+\|\sqrt{r}\p_z v_r\|_2^2)\nn \\
&&\qquad +\nu_2(\f{1}{r}\p_r(rv_r)+\p_z v_z,\p_r(r v_r))
=(g_1, \bar\rho r v_r)+(\f{2M_0}{r}v_\theta, \bar\rho v_r),
\label{cc2.1}\\
&&\f{1}{2}\f{d}{dt}\|\sqrt{\bar\rho r}v_\theta\|_2^2
+\nu_1(\|\f{v_\theta}{\sqrt{r}}\|_2^2+\|\sqrt{r}\p_r v_\theta\|_2^2+\|\sqrt{r}\p_z v_\theta\|_2^2)
=(g_2, \bar\rho r v_\theta)\label{cc2.2}
\eea
and
\bea
&&\f{1}{2}\f{d}{dt}\|\sqrt{\bar\rho r}v_z\|_2^2
+\gamma (\p_z({\bar\rho}^{\gamma-2} \phi), \bar\rho r v_z)
+\nu_1(\|\sqrt{r}\p_r v_z\|_2^2+\|\sqrt{r}\p_z v_r\|_2^2)\nn \\
&&\qquad +\nu_2(\f{1}{r}\p_r(rv_r)+\p_z v_z, r\p_z  v_r)
=(g_3, \bar\rho r v_z).
\label{cc2.3}
\eea
Adding \eqref{cc2.1} and \eqref{cc2.3} yields
\bea
&&\f{1}{2}\f{d}{dt}(\|\sqrt{{\bar\rho} r}v_r\|_2^2+\|\sqrt{{\bar\rho} r}v_z\|_2^2)
-\gamma ({\bar\rho}^{\gamma-2} \phi, \p_r({\bar\rho}r v_r)+\p_z({\bar\rho} r v_z))\nn\\
&&\quad +\nu_1(\|\f{v_r}{\sqrt{r}}\|_2^2+\|\sqrt{r}\p_r v_r\|_2^2+\|\sqrt{r}\p_z v_r\|_2^2
+\|\sqrt{r}\p_r v_z\|_2^2+\|\sqrt{r}\p_z v_z\|_2^2)+\nu_2\|\sqrt{r}(\f{1}{r}\p_r(rv_r)+\p_z v_z)\|_2^2\nn \\
&&\quad
=(g_1, {\bar\rho} r v_r)+(\f{2M_0}{r}v_\theta, {\bar\rho}v_r)+(g_3, {\bar\rho} r v_z).\label{cc2.4}
\eea
By $\ds|(\f{2M_0}{r}v_\theta, {\bar\rho}v_r)|\le \f{\nu_1}{2}\|\f{v_r}{\sqrt{r}}\|_2^2  +C\|\f{v_\theta}{\sqrt{r}}\|_2^2
$, one has that from \eqref{cc2.2} and \eqref{cc2.4},
\bea
&&\f{1}{2}\f{d}{dt}\|\sqrt{{\bar\rho} r}v\|_2^2
-\gamma ({\bar\rho}^{\gamma-2} \phi, \p_r({\bar\rho} r v_r)+\p_z({\bar\rho} r v_z))\nn\\
&&\quad +\nu_1(\|\f{v_r}{\sqrt{r}}\|_2^2+\|\f{v_\theta}{\sqrt{r}}\|_2^2+\|\sqrt{r}\p_r v\|_2^2+\|\sqrt{r}\p_z v\|_2^2)+\nu_2\|\sqrt{r}(\f{1}{r}\p_r(rv_r)+\p_z v_z)\|_2^2
\nn \\
&&
\le C(g_1, {\bar\rho} r v_r)+C(g_2, {\bar\rho} r v_\theta)+C(g_3, {\bar\rho} r v_z).\label{cc2.5}
\eea
In addition, it follows from $\ds\int_\Omega (1.13)\times \gamma \bar \rho^{\gamma-2} r \phi drdz$ that
\be
\f{\gamma}{2}\f{d}{dt}\|\sqrt{{\bar\rho}^{\gamma-2} r}\phi\|_2^2 +\gamma ({\bar\rho}^{\gamma-2} \phi, \p_r({\bar\rho}r v_r)
+\p_z({\bar\rho} r v_z))=\gamma (f, {\bar\rho}^{\gamma-2} r \phi).\label{cc2.6}
\ee
Consequently, adding \eqref{cc2.5}, \eqref{cc2.6} and integrating with respect to the time variable $\tau$ over $(0,t)$  yields Lemma 2.1. \ef

\begin{lemma}\label{pcl2} ({\bf Weighted $L^2-$estimate of $(\p_t\phi, \p_tv, Dv)$}). For the solution
$(\phi, v)\in C([0,\infty),\tilde H^2\times\tilde H^3)$ of problem (1.13)-(1.16) with (1.17)-(1.18),
we have
\bea
&& \|\sqrt{r}D v \|_2^2 + \|\sqrt{r}(\f{1}{r}\p_r(r v_r)+\p_z v_z)\|_2^2   + \int_0^t\biggl( \|\sqrt{{\bar\rho} r }\p_t v\|_2^2+ \|\sqrt{ r {\bar\rho}^{\gamma-2}}{\p_t\phi}\|_2^2\biggl)d\tau \nn\\
&&\le C \|(\phi_0,v_0)\|_{\tilde H^1}^2+ C \int_0^t (A_1+A_2)d\tau,\nn
\eea
where $A_1$ is defined in Lemma 2.1, and $A_2=|(g, {\bar\rho}r\p_t v )|+|(f, r{\bar\rho}^{\gamma-2}\p_t\phi)|$.
\end{lemma}
\pf By computing $\ds\int_\Omega (1.14)\times\bar\rho r \p_t v_rdrdz$, $\ds\int_\Omega\eqref{cc1.7}\times \bar\rho r \p_t v_\theta drdz$  and $\ds\int_\Omega \eqref{cc1.8}\times\bar\rho r \p_t v_zdrdz$, we obtain that
\bea
&& (\p_t v_r, {\bar\rho}r \p_t v_r)-(\f{2M_0}{r}v_\theta, {\bar\rho}\p_t v_r) +\gamma (\p_r ({\bar\rho}^{\gamma-2}\phi),
{\bar\rho}r \p_t v_r) \nn\\
&&\qquad -\nu_1( \p_r (\f{1}{r}\p_r(rv_r))+\p_z^2 v_r,  r \p_t v_r)-\nu_2(\p_r(\f{1}{r}\p_r(rv_r)+\p_z v_z), r\p_t v_r)
\nn\\
&&\quad =(g_1, {\bar\rho}r \p_t v_r),\label{cc2.8}\\
&& (\p_t v_\theta, {\bar\rho}r \p_t v_\theta)
 -\nu_1(\p_r (\f{1}{r}\p_r(rv_\theta))+\p_z^2 v_\theta,  r \p_t v_\theta) =(g_2, {\bar\rho}r \p_t v_\theta),\label{cc2.9}
\eea
and
\bea
&& (\p_t v_z, {\bar\rho}r \p_t v_z)+\gamma (\p_z ({\bar\rho}^{\gamma-2}\phi), {\bar\rho}r \p_t v_z) \nn\\
&&\quad -\nu_1(\p_r^2 v_z +\f{1}{r}\p_r v_z+\p_z^2 v_z,  r \p_t v_z)-\nu_2(\p_z(\f{1}{r}\p_r(rv_r)+\p_z v_z), r\p_t v_z)
\nn\\
&&=(g_3, {\bar\rho}r \p_t v_z).\label{cc2.10}
\eea
Note that
\bea
&&-(\f{\nu_1}{{\bar\rho}}\p_r (\f{1}{r}\p_r(rv_r)),{\bar\rho} r \p_t v_r)
=\f{\nu_1}{2}\f{d}{dt}(\|\f{v_r}{r^{1/2}}\|_2^2+\|r^{1/2}\p_r v_r\|_2^2).\nn
\eea
This, together with \eqref{cc2.8}-\eqref{cc2.10}, derives
\bea
&& \|\sqrt{{\bar\rho} r }\p_t v\|_2^2-\gamma ({\bar\rho}^{\gamma-2}\phi, \p_r({\bar\rho} r \p_t v_r)+\p_z({\bar\rho} r \p_t v_z)) \nn\\
&&\quad +\f{d}{dt}\biggl(\f{\nu_1}{2}(\|\f{v_r }{\sqrt{r}}\|_2^2+\|\f{v_\theta }{\sqrt{r}}\|_2^2+\|\sqrt{r}D v \|_2^2)+\f{\nu_2}{2}\|\sqrt{r}(\f{1}{r}\p_r(r v_r)+\p_z v_z)\|_2^2 \biggl)\nn\\
&&\le C(g, {\bar\rho} r \p_t v )+C(\f{2M_0}{r}v_\theta,{\bar\rho}   \p_t v_r).\label{cc2.11}
\eea
Computing $\ds\int_\Omega (1.13)\times\gamma{\bar\rho}^{\gamma-2} r\p_t\phi drdz$ yields
\be
\gamma\|\sqrt{ r {\bar\rho}^{\gamma-2}}{\p_t\phi}\|_2^2 +\gamma(\f{1}{r}\p_r(r {\bar\rho} v_r)+\p_z (\bar \rho v_z), r {\bar\rho}^{\gamma-2}{\p_t\phi})=
\gamma(f, r {\bar\rho}^{\gamma-2} {\p_t\phi}). \label{cc2.12}
\ee
In addition, we see that
\bea
&&\gamma(\f{1}{r}\p_r(r {\bar\rho} v_r)+\p_z (\bar \rho v_z), r {\bar\rho}^{\gamma-2}{\p_t\phi})
-\gamma(\p_t\big(\f{1}{r}\p_r(r {\bar\rho} v_r)+\p_z (\bar \rho v_z)\big), r {\bar\rho}^{\gamma-2}\phi )\nn\\
&&=-\gamma \p_t( {\bar\rho}^{\gamma-2}\phi, \p_r(r{\bar\rho} v_r)+r\p_z (\bar \rho v_z))
+2\gamma ( {\bar\rho}^{\gamma-2}{\p_t\phi}, \p_r(r{\bar\rho} v_r)+r\p_z (\bar \rho v_z)).\label{cc2.13}
\eea
Then adding \eqref{cc2.11} and \eqref{cc2.12},  and using \eqref{cc2.13}, we arrive at
\bea
&& \|\sqrt{{\bar\rho} r }\p_t v\|_2^2+\gamma\|\sqrt{ r {\bar\rho}^{\gamma-2}}{\p_t\phi}\|_2^2
+\f{d}{dt}\biggl(\f{\nu_1}{2}(\|\f{v_r }{\sqrt{r}}\|_2^2+\|\f{v_\theta }{\sqrt{r}}\|_2^2+\|\sqrt{r}D v \|_2^2)\nn\\
&&\quad +\f{\nu_2}{2}\|\sqrt{r}(\f{1}{r}\p_r(r v_r)+\p_z v_z)\|_2^2
-\gamma  ( {\bar\rho}^{\gamma-2}\phi, \p_r(r{\bar\rho} v_r)+r\p_z (\bar \rho v_z)) \biggl)\nn\\
&&\le (g ,{\bar\rho} r \p_t v ) +(\f{2M_0v_\theta}{r},{\bar\rho} \p_t v_r)+\gamma(f, r {\bar\rho}^{\gamma-2} {\p_t\phi})-2\gamma ( {\bar\rho}^{\gamma-2}{\p_t\phi}, \p_r(r{\bar\rho} v_r)+r\p_z (\bar \rho v_z)).\nn\\
&&\label{cc2.14}
\eea
Since
\bea
&&\gamma  ( {\bar\rho}^{\gamma-2}\p_t\phi, \p_r(r{\bar\rho} v_r)+r\p_z (\bar \rho v_z))\nn\\
&&\quad \le C\|\sqrt{  r }\biggl(\f{1}{r} \p_r(r  v_r)+\p_z v_z\biggl)\|_2^2 +C(\|r^{1/2}\p_r v_r \|_2^2+ \|\f{v_r}{\sqrt{r}}\|_2^2)+\f{\gamma}{4} \|\sqrt{ r {\bar\rho}^{\gamma-2}}{\p_t\phi}\|_2^2,
\label{cc2.15}
\eea
and
\bea
&&|({\bar\rho}^{\gamma-2}\phi, \p_r(r{\bar\rho} v_r)+r\p_z(\bar \rho  v_z))|\nn\\
&&\le \f{\nu_2}{4}\|\sqrt{  r }\biggl(\f{1}{r}(\p_r(r  v_r)+\p_z v_z\biggl)\|_2^2 +\f{\nu_1}{4}\|r^{1/2}\p_r v_r \|_2^2+C \|\sqrt{ r {\bar\rho}^{\gamma-2}}\phi\|_2^2+C\|\f{v_r}{\sqrt{r}}\|_2^2,
\label{cc2.16}
\eea
integrating \eqref{cc2.14} with respect to the time variable $\tau$ over $(0,t)$  and combining \eqref{cc2.15}-\eqref{cc2.16}  yield
\bea
&& \f{1}{2}\int_0^t \big(\|\sqrt{{\bar\rho} r }\p_t v\|_2^2+\f{\gamma}{2}\|\sqrt{ r {\bar\rho}^{\gamma-2}}{\p_t\phi}\|_2^2\big)d\tau
+ \biggl(\f{\nu_1}{4}(\|\f{v_r }{\sqrt{r}}\|_2^2+\|\f{v_\theta }{\sqrt{r}}\|_2^2+\|\sqrt{r}D v \|_2^2) \nn\\
&&\quad +\f{\nu_2}{4}\|\sqrt{r}(\f{1}{r}\p_r(r v_r)+\p_z v_z)\|_2^2-C \|\sqrt{ r {\bar\rho}^{\gamma-2}}\phi\|_2^2 -C\|\f{v_r}{\sqrt{r}}\|_2^2\biggl)\nn\\
&&\le  C \|(\phi_0,v_0)\|_{\tilde H^1}^2+\int_0^t \biggl(|(g ,{\bar\rho} r \p_t v )| +|(\f{2M_0v_\theta}{r},{\bar\rho} \p_t v_r)|+|\gamma(f, r {\bar\rho}^{\gamma-2} {\p_t\phi})|+\nn\\
&&\quad C\|\sqrt{  r }\big(\f{1}{r} \p_r(r  v_r)+\p_z v_z\big)\|_2^2 +C\|r^{1/2}\p_r v_r \|_2^2+C\|\f{v_r}{\sqrt{r}}\|_2^2 \biggl)d\tau.\label{cc2.17}
\eea
Together with Holder's inequality and Lemma \ref{pcl1}, this yields
\bea
&& \|\sqrt{r}D v \|_2^2 + \|\sqrt{r}(\f{1}{r}\p_r(r v_r)+\p_z v_z)\|_2^2   + \int_0^t\biggl( \|\sqrt{{\bar\rho} r }\p_t v\|_2^2+ \|\sqrt{ r {\bar\rho}^{\gamma-2}}{\p_t\phi}\|_2^2\biggl)d\tau \nn\\
&&\le C\|(\phi_0,v_0)\|_{\tilde H^1}^2+ C \int_0^t (A_1+A_2)d\tau,\label{cc2.18}
\eea
which completes the proof of Lemma \ref{pcl2}.\ef

Taking $\p_t^k\p_z^j$ ($k=0,1$ and $j=0,1,2$) on both sides of equations (1.13)-\eqref{cc1.8}, we then have
\bea
&&\p_t (\p_t^k\p_z^j\phi) +\f{1}{r}\p_r (r {\bar\rho} \p_t^k\p_z^jv_r)+\bar \rho \p_z\p_t^k\p_z^j v_z=\p_t^k\p_z^jf, \label{cc2.19}\\
&&\p_t(\p_t^k\p_z^j v_r) -\f{2M_0}{r^2}\p_t^k\p_z^jv_\theta +\gamma \p_r( {\bar\rho}^{\gamma-2}  \p_t^k\p_z^j\phi)
-\f{\nu_1}{{\bar\rho}}\bigg(\p_r(\f{1}{r}\p_r(r\p_t^k\p_z^jv_r))+\p_z^2 \p_t^k\p_z^jv_r\biggl)\nn\\
&&\quad \quad \quad -\f{\nu_2}{{\bar\rho}}\p_r(\f{1}{r}\p_r (r  \p_t^k\p_z^jv_r)+\p_z  \p_t^k\p_z^jv_z )=\p_t^k\p_z^jg_1,
 \label{cc2.20}\\
&&\p_t(\p_t^k\p_z^j v_\theta)
-\f{\nu_1}{{\bar\rho}}\biggl(\p_r(\f{1}{r}\p_r(r\p_t^k\p_z^jv_\theta))+\p_z^2 \p_t^k\p_z^jv_\theta\biggl)=\p_t^k\p_z^jg_2,\label{cc2.21} \\
&&\p_t(\p_t^k\p_z^j v_z)+\gamma \p_z ({\bar\rho}^{\gamma-2} \p_t^k\p_z^j\phi)-\f{\nu_1}{{\bar\rho}}\biggl( \p_r^2 \p_t^k\p_z^jv_z +\p_z^2\p_t^k\p_z^jv_z +\f{1}{r}\p_r \p_t^k\p_z^jv_z \biggl)\nn\\
&&\quad\quad\quad-\f{\nu_2}{{\bar\rho}}\p_z(\f{1}{r}\p_r (r  \p_t^k\p_z^jv_r)+\p_z \p_t^k\p_z^jv_z)=\p_t^k\p_z^jg_3.\label{cc2.22}
\eea
Adding $ \ds\int_\Omega \gamma {\bar\rho}^{\gamma-2} r\p_t^k\p_z^j \phi \times \eqref{cc2.19}drdz$,$\ds\int_\Omega {\bar\rho} r\p_t^k\p_z^j v_r\times \eqref{cc2.20}drdz$,$\ds\int_\Omega  {\bar\rho} r\p_t^k\p_z^j v_\theta \times \eqref{cc2.21}drdz$ and
\newline $\ds\int_\Omega    {\bar\rho} r\p_t^k\p_z^j v_z \times \eqref{cc2.22}drdz$, then as in Lemma 2.1-Lemma 2.2, we can obtain

\begin{lemma}\label{pcl3} ({\bf Weighted $L^2-$estimate of $(\p_t^k\p_z^j\phi, \p_t^k\p_z^jv)$}
with $k=0,1$ and $j=0,1,2$). For the solution
$(\phi, v)\in C([0,\infty),\tilde H^2\times\tilde H^3)$ of problem (1.13)-(1.16) with (1.17)-(1.18),
we have
\bea
&&  \|\sqrt{{\bar\rho} r}\p_t^k\p_z^jv \|_2^2 +\|\sqrt{{\bar\rho}^{\gamma-2} r}\p_t^k\p_z^j\phi\|_2^2  + \int_0^t \biggl( \|\f{\p_t^k\p_z^jv_r}{\sqrt{r}}\|_2^2+\|\f{\p_t^k\p_z^jv_\theta}{\sqrt{r}}\|_2^2\nn\\
&&\quad +\|\sqrt{r}\p_r\p_t^k\p_z^j v\|_2^2+\|\sqrt{r}\p_t^k\p_z^{j+1}v\|_2^2+\|\sqrt{r}(\f{1}{r}\p_r(r\p_t^k\p_z^jv_r)+ \p_t^k\p_z^{j+1}v_z)\|_2^2\biggl) d\tau\nn\\
&&\le C\|\sqrt{r}
\p_t^{k_1}\p_z^{j_1}v_0\|_2^2+\|\sqrt{r}
\p_t^{k_1}\p_z^{j_1}\phi_0\|_2^2+ C\sum\limits_{k_1=0,j_1=0}^{k,j}\biggl(\int_0^t
|(\p_t^{k_1}\p_z^{j_1}g,{\bar\rho} r \p_t^{k_1}\p_z^{j_1}v)|\nn\\
&&\quad
+|(\p_t^{k_1}\p_z^{j_1}f,{\bar\rho}^{\gamma-2} r\p_t^{k_1}\p_z^{j_1}\phi)|\biggl)d\tau.\nn
\eea
\end{lemma}

On the other hand, by  equations \eqref{cc2.19}-\eqref{cc2.22} and Lemma \ref{pcl3} with $k=1,j=0$,
as in Lemma \ref{pcl2},  we can also obtain

\begin{lemma}\label{pcl4} ({\bf Weighted $L^2-$estimate of $(\p_t^2\phi, \p_t^2v, \p_tDv)$}). For the solution
$(\phi, v)\in C([0,\infty),\tilde H^2\times\tilde H^3)$ of problem (1.13)-(1.16) with (1.17)-(1.18),
we have
\bea
&&   \|\sqrt{r}{\p_{ t} Dv} \|_2^2 + \|\sqrt{r}(\f{1}{r}\p_r(r \p_tv_r)+\p_{zt}^2v_z)\|_2^2 +\int_0^t\big( \|\sqrt{{\bar\rho} r }\p_t^2 v\|_2^2+ \|\sqrt{ r {\bar\rho}^{\gamma-2}}\p_t^2\phi\|_2^2\big)d\tau \nn\\
&&\le C (\| \phi_0\|_{\tilde H^2}^2+\|v_0\|_{\tilde H^1}^3)+C\int_0^t (A_1+A_2+A_3)d\tau,\nn
\eea
where $A_1$ and $A_2$ are defined in Lemma 2.1 and Lemma 2.2 respectively, and
$A_3=|(\p_tg,{\bar\rho} r \p_t v )|+|(\p_tf, r  {\bar\rho}^{\gamma-2} \p_t\phi)|+|(\p_tg, {\bar\rho} r \p_t^2 v )|
+ |(\p_tf, r  {\bar\rho}^{\gamma-2} {\p_t^2\phi})|$.
\end{lemma}

Next we start to derive the a priori estimates of $\phi$.
By $(\nu_1+\nu_2)\p_r\bigl({\bar\rho}^{\gamma}\times (1.13)\bigl) $ and direct computation, one arrives at
\bea
&&(\nu_1+\nu_2)\p_{rt}^2({\bar\rho}^{\gamma-2}\phi)+(\nu_1+\nu_2) {\bar\rho}^{\gamma-1}\p_r(\f{1}{r}\p_r(r v_r))\nn\\ &&=(\nu_1+\nu_2)\bigl\{\p_r({\bar\rho}^{\gamma-2}(f -v_r\p_r{\bar\rho})) -\p_r({\bar\rho}^{\gamma-1})(\f{1}{r}\p_r(r v_r)+\p_zv_z)-{\bar\rho}^{\gamma-1}\p_{rz}^2v_z \bigl\}.\nn\\
&&\label{cc2.23}
\eea
Adding \eqref{cc2.23} and ${\bar\rho}^{\gamma-2}\times (1.13)$ yields
\be
(\nu_1+\nu_2)\p_{tr}^2({\bar\rho}^{\gamma-2}\phi)+\gamma {\bar\rho}^{\gamma}\p_r({\bar\rho}^{\gamma-2}\phi)=h, \label{cc2.25}
\ee
where
\bea
h& =&(\nu_1+\nu_2)\biggl\{\p_r({\bar\rho}^{\gamma-2}(f -v_r\p_r{\bar\rho})) -\p_r({\bar\rho}^{\gamma-1})(\f{1}{r}\p_r(r v_r)+\p_zv_z) \biggl\}\nn\\
&& +{\bar\rho}^{\gamma}(-\p_t v_r +\f{2M_0}{r^2}v_\theta)
+{\bar\rho}^{\gamma-1}(\nu_1\p_z^2v_r-\nu_1\p_{rz}^2 v_z)+{\bar\rho}^{\gamma}g_1
\nn\\
&=& h_1+h_2\nn
\eea
with $h_1=(\nu_1+\nu_2)\p_r({\bar\rho}^{\gamma-2}f)$ and $h_2=h-h_1$.

By $\ds\int_\Omega r\p_r({\bar\rho}^{\gamma-2}\phi)\times\eqref{cc2.25} drdz$, we have
\be
\f{\nu_1+\nu_2}{2}\f{d}{dt}\|r^{1/2}\p_r({\bar\rho}^{\gamma-2}\phi)\|_2^2+\|\gamma\sqrt{{\bar\rho}^{\gamma} r}\p_r({\bar\rho}^{\gamma-2}\phi)\|_2^2 \le (\p_r({\bar\rho}^{\gamma-2}\phi),rh).\label{cc2.26}
\ee
Also, we obtain that from  $\ds\int_\Omega  r\p_r^2({\bar\rho}^{\gamma-2}\phi)\times \p_r\eqref{cc2.25}drdz$,

\bea
&&\f{\nu_1+\nu_2}{2}\f{d}{dt}\|r^{1/2}\p_r^2({\bar\rho}^{\gamma-2}\phi)\|_2^2+\|\sqrt{{\bar\rho}^{\gamma} r}\p_r^2({\bar\rho}^{\gamma-2}\phi)\|_2^2 \nn\\
&&\le (\p_r^2({\bar\rho}^{\gamma-2}\phi),r\p_rh )
+(\p_r^2({\bar\rho}^{\gamma-2}\phi),r\p_r({\bar\rho}^{\gamma} )\p_r  ({\bar\rho}^{\gamma-2}\phi)).
\label{cc2.28}
\eea
Next, we show that

\begin{lemma}\label{pcl5}  ({\bf Weighted higher order energy estimate}). For the solution
$(\phi, v)\in C([0,\infty),\tilde H^2\times\tilde H^3)$ of problem (1.13)-(1.16) with (1.17)-(1.18),
we have
\bea
&& \|r^{1/2}\p_r({\bar\rho}^{\gamma-2}\phi)\|_2^2+\|r^{1/2}\p_r({\bar\rho}^{\gamma-2}\p_z\phi)\|_2^2+\|r^{1/2}\p_r^2({\bar\rho}^{\gamma-2}\phi)\|_2^2 +\int_0^t \biggl(\|r^{1/2}\p_r({\bar\rho}^{\gamma-2}\phi)\|_2^2\nn\\
&&\quad +\|r^{1/2} \p_r^2v_r\|_2^2+\|r^{1/2}\p_r({\bar\rho}^{\gamma-2}\p_z\phi)\|_2^2
+\|r^{1/2}\p_r^2\p_z v_r\|_2^2+\|\sqrt{  r}\p_r^2({\bar\rho}^{\gamma-2}\phi)\|_2^2 \biggl)d\tau\nn\\
&&\le  C(\|\phi_0\|_{\tilde H^2}^2+\|v_0\|_{\tilde H^3}^2)+C\int_0^t \biggl(A_1+A_2+\|  g\|_{\tilde H^1}^2+\|f\|_{\tilde H^1}^2 +|(r\p_r({\bar\rho}^{\gamma-2}\p_z\phi ), {\bar\rho}^{\gamma-2}\p_{rz}^2f)|\nn\\
 &&\quad
 +|(\p_r^2({\bar\rho}^{\gamma-2}\phi), r{\bar\rho}^{\gamma-2}\p_r^2f)|
 +|(\p_z^{2}g,{\bar\rho} r \p_z^{2}v)|
 +|(\p_z^{2}f,{\bar\rho}^{\gamma-2} r \p_z^{2}\phi)|\biggl)d\tau.
\nn
\eea
\end{lemma}
\pf We see that
\bea
&&(\p_r({\bar\rho}^{\gamma-2}\phi),rh)\le \f{1}{2}\|\sqrt{{\bar\rho}^{\gamma} r}\p_r({\bar\rho}^{\gamma-2}\phi)\|_2^2 +C\|r^{1/2}h\|_2^2. \label{cc2.29}
\eea
In addition,
\bea
\|r^{1/2}h\|_2^2 & \le& C\biggl\{ \|r^{1/2}f\|_2^2+
\|r^{-1/2}v_r\|_2^2 +\|r^{1/2}\p_r v_r\|_2^2+\|r^{1/2}\p_z v_z\|_2^2 \nn\\
&&  +\|r^{1/2}\p_{rz}^2 v_z\|_2^2+\|r^{1/2}\p_z^2 v_z\|_2^2+\|r^{1/2}\p_t v_r \|_2^2+\|r^{1/2}g_1\|_2^2+ \|r^{1/2}\p_rf\|_2^2
\biggl\}.\nn\\
& \le & C\biggl\{ \|f\|_{\tilde H^1}^2+\|g\|_{\tilde H^1}^2 +\|r^{-1/2}v_r\|_2^2 +\|r^{1/2}Dv\|_2^2+\|r^{1/2}\p_{z}  Dv_z\|_2^2+\|r^{1/2}\p_t v_r \|_2^2\biggl\}.\nn
\eea
Together with  \eqref{cc2.26}, this yields
\bea
&&\f{\nu_1+\nu_2}{2}\f{d}{dt}\|r^{1/2}\p_r({\bar\rho}^{\gamma-2}\phi)\|_2^2+\f{\gamma}{2}\|\sqrt{{\bar\rho}^{\gamma} r}\p_r({\bar\rho}^{\gamma-2}\phi)\|_2^2 \nn \\
&&\le C\biggl\{ \|f\|_{\tilde H^1}^2+\|g\|_{\tilde H^1}^2 +\|r^{-1/2}v_r\|_2^2 +\|r^{1/2}Dv\|_2^2+\|r^{1/2}\p_{z}  Dv_z\|_2^2+\|r^{1/2}\p_t v_r \|_2^2\biggl\}. \label{cc2.30}
\eea
Note that
\bea
&&(\p_r^2({\bar\rho}^{\gamma-2}\phi),r\p_rh )\le
\f{1}{8}\|\sqrt{{\bar\rho}^{\gamma} r}\p_r^2({\bar\rho}^{\gamma-2}\phi)\|_2^2+C\|r^{1/2}\p_r h_2\|_2^2
+(\p_r^2({\bar\rho}^{\gamma-2}\phi),r\p_rh_1),\nn\\
&&(\p_r^2({\bar\rho}^{\gamma-2}\phi),r\p_r({\bar\rho}^{\gamma} )\p_r  ({\bar\rho}^{\gamma-2}\phi))
\le \f{1}{4}\|\sqrt{{\bar\rho}^{\gamma} r}\p_r^2({\bar\rho}^{\gamma-2}\phi)\|_2^2+C \|r^{1/2}\p_r({\bar\rho}^{\gamma} )\p_r  ({\bar\rho}^{\gamma-2}\phi)\|_2^2,\nn
\eea
moreover,
\bea
&&\|r^{1/2}\p_r h_2\|_2^2\nn\\
&&\le
C\biggl(
\|r^{-1/2}v_r\|_2^2+\|r^{-1/2}v_\theta\|_2^2 +\|r^{1/2}D v \|_2^2+\|r^{1/2}\p_t v_r \|_2^2+
\|r^{1/2}\p_{rz}^2v_z\|_2^2+\|r^{1/2}\p_{z}^2v_r\|_2^2\nn\\
&&\quad +\|r^{1/2}\p_r^2 v_r\|_2^2 +\|r^{1/2}\p_{rt}^2v_r \|_2^2+\|r^{1/2}\p_z^2\p_rv_r \|_2^2+\|r^{1/2}\p_r^2 \p_zv_z\|_2^2+\|g\|_{\tilde H^1}^2\biggl)\label{cc2.31}
\eea
and
\bea
&&(\p_r^2({\bar\rho}^{\gamma-2}\phi),r\p_rh_1 )\nn\\
&&=(\nu_1+\nu_2)(\p_r^2({\bar\rho}^{\gamma-2}\phi), r\p_r^2({\bar\rho}^{\gamma-2})f
+r\p_r{\bar\rho}^{\gamma-2}\p_rf
+r{\bar\rho}^{\gamma-2}\p_r^2f)\nn\\
&&\le \f{1}{8}\|\sqrt{{\bar\rho}^{\gamma} r}\p_r^2({\bar\rho}^{\gamma-2}\phi)\|_2^2+C\| f\|_{\tilde H^1}^2 +
C(\p_r^2({\bar\rho}^{\gamma-2}\phi), r{\bar\rho}^{\gamma-2}\p_r^2f).\nn
\eea
Together with \eqref{cc2.28}, this yields
\bea
&&\f{d}{dt}\|r^{1/2}\p_r^2({\bar\rho}^{\gamma-2}\phi)\|_2^2+\|\sqrt{{\bar\rho}^{\gamma} r}\p_r^2({\bar\rho}^{\gamma-2}\phi)\|_2^2 \nn\\
&&\le
C\biggl(
\|r^{-1/2}v_r\|_2^2+\|r^{-1/2}v_\theta\|_2^2 +\|r^{1/2}D v \|_2^2+\|r^{1/2}\p_t v_r \|_2^2+
\|r^{1/2}\p_{rz}^2v_z\|_2^2+\|r^{1/2}\p_{z}^2v_r\|_2^2\nn\\
&&\quad +\|r^{1/2}\p_r^2 v_r\|_2^2 +\|r^{1/2}\p_{rt}^2v_r \|_2^2+\|r^{1/2}\p_z^2\p_rv_r \|_2^2+\|r^{1/2}\p_r^2 \p_zv_z\|_2^2+\|g\|_{\tilde H^1}^2+ \| f\|_{\tilde H^1}^2  \nn\\
&&\quad+
 (\p_r^2({\bar\rho}^{\gamma-2}\phi), r{\bar\rho}^{\gamma-2}\p_r^2f)\biggl).\label{cc2.32}
\eea
On the other hand, we rewrite \eqref{cc2.20} with $k=0$ as
\bea
(\nu_1+\nu_2)\p_z^j\p_r^2 v_r&={\bar\rho}\p_t\p_z^j v_r -{\bar\rho}\f{2M_0}{r^2}\p_z^jv_\theta +\gamma {\bar\rho}\p_r({\bar\rho}^{\gamma-2} \p_z^j\phi)
- (\nu_1+\nu_2) \p_r(\f{\p_z^jv_r}{r} )\nn\\
&-\nu_1 \p_z^{2+j}v_r-\nu_2\p_r\p_z^{j+1}v_z-{\bar\rho}\p_z^j g_1.\label{cc2.33}
\eea
This derives that
\bea
&&(\nu_1+\nu_2)\|r^{1/2} \p_r^2\p_z^j v_r\|_2^2\nn\\
&& \le C\biggl\{\|r^{1/2}\p_t\p_z^j v_r\|_2^2 +\|\f{1}{r^{3/2}}\p_z^jv_\theta\|_2^2
+\|r^{1/2}\p_r({\bar\rho}^{\gamma-2} \p_z^j\phi)\|_2^2
+ \| r^{1/2}\p_r(\f{\p_z^jv_r}{r} )\|_2^2\nn\\
&&\quad +\|r^{1/2}\p_z^{2+j}v_r\|_2^2+\|r^{1/2}\p_r\p_z^{j+1}v_z\|_2^2+\|r^{1/2}\p_z^j g_1\|_2^2\biggl\}.\label{cc2.34}
\eea
Combining \eqref{cc2.30}, \eqref{cc2.34} with $j=0$, Lemma \ref{pcl1}-\ref{pcl2} and Lemma \ref{pcl3} with $k=0,j=1$  and integrating with respect to the time variable $\tau$ over$(0,t)$,  we have
\bea
&& \|r^{1/2}\p_r({\bar\rho}^{\gamma-2}\phi)\|_2^2+\int_0^t \big(\|r^{1/2}\p_r({\bar\rho}^{\gamma-2}\phi)\|_2^2+\|r^{1/2} \p_r^2v_r\|_2^2\big)d\tau\nn\\
&&\le C \|(\phi_0,v_0)\|_{\tilde H^2}^2  +C\int_0^t \biggl(  A_1+A_2+\|  g\|_{\tilde H^1}^2+\| f\|_{\tilde H^1}^2\biggl)d\tau.\nn\\
 &&
 \label{cc2.35}
\eea
Note that we have from \eqref{cc2.25}
\be
(\nu_1+\nu_2)\p_{tr}^2({\bar\rho}^{\gamma-2}\p_z\phi)+{\bar\rho}^{\gamma}\gamma \p_r({\bar\rho}^{\gamma-2}\p_z\phi) =\p_zh. \label{cc2.36}
\ee
Then it follows from $
\ds\int_\Omega\p_r({\bar\rho}^{\gamma-2}\p_z\phi)\times \eqref{cc2.36} rdrdz$ that
\bea
&&\f{d}{dt}\|r^{1/2}\p_r({\bar\rho}^{\gamma-2}\p_z\phi)\|_2^2+\|r^{1/2}\p_r({\bar\rho}^{\gamma-2}\p_z\phi)\|_2^2\nn\\
&&\le  C(\p_zh,r\p_r({\bar\rho}^{\gamma-2}\p_z\phi ))\nn\\
&&\le C \biggl( \|r^{1/2}\p_z g\|_2^2+\|r^{1/2}\p_rf\|_2^2+({\bar\rho}^{\gamma-2}r\p_r({\bar\rho}^{\gamma-2}\p_z\phi ),{\bar\rho}^{\gamma-2}\p_{rz}^2f)+\|r^{1/2}\p_{rz}^2v_r\|_2^2\nn\\
&&\quad+\|r^{1/2}\p_{z}v_r\|_2^2
+\|r^{1/2}\p_{r}v_r\|_2^2++\|r^{1/2}\p_{z}^2v_z\|_2^2
+\|r^{1/2}\p_{tz}^2v_r\|_2^2+\|r^{1/2}\p_{z}v_\theta\|_2^2\nn\\
&&\quad+\|r^{1/2}\p_z^3v_z\|_2^2
\biggl).
 \label{cc2.37}
\eea
This, together with \eqref{cc2.34} with $j=1$, yields
\bea
&& \|r^{1/2}\p_r({\bar\rho}^{\gamma-2}\p_z\phi)\|_2^2+\int_0^t \big(\|r^{1/2}\p_r({\bar\rho}^{\gamma-2}\p_z\phi)\|_2^2+\|r^{1/2}\p_r^2\p_z v_r\|_2^2\big)d\tau\nn\\
&&\le C(\|\phi_0\|_{\tilde H^2}^2+\|v_0\|_{\tilde H^3}^2)+\int_0^t C \biggl( (r\p_r({\bar\rho}^{\gamma-2}\p_z\phi ),{\bar\rho}^{\gamma-2}\p_{rz}^2f)+\|r^{1/2}\p_z g\|_2^2+\|r^{1/2}\p_rf\|_2^2\nn\\
&&\quad +\|r^{1/2}\p_{rz}^2v_r\|_2^2
+\|r^{1/2}D v\|_2^2+\| \p_{z}^2v_z\|_{\tilde H^1}^2
+\|r^{1/2}\p_{tz}^2v_r\|_2^2+\|r^{1/2}\p_z^{3}v_r\|_2^2
\biggl) d\tau. \label{cc2.38}
\eea
Consequently, by \eqref{cc2.32}, \eqref{cc2.35}, \eqref{cc2.38}, Lemma \ref{pcl1}-\ref{pcl2}, and Lemma \ref{pcl3} with $k=0,j=1,2$, we
complete the proof of Lemma 2.5.\ef

\section{Global estimates and proof of Theorem 1.1}
For the solution $(\phi, v)\in C([0,\infty),\tilde H^2\times\tilde H^3)$ of problem (1.13)-(1.16) with (1.17)-(1.18),
by Lemma \ref{pcl1}-\ref{pcl5}, we obtain that
\bea
&&
\|v\|_{\tilde H^1}^2+\|\p_t v\|_{\tilde H^1}^2+\|\sqrt{r}\phi\|_2^2+\|\sqrt{r}\p_t \phi\|_2^2
+
\|r^{1/2}\p_r({\bar\rho}^{\gamma-2}\phi)\|_2^2+\|r^{1/2}\p_rD({\bar\rho}^{\gamma-2} \phi)\|_2^2
 \nn\\
&&\quad +\int_0^t \biggl(\|\f{v_r}{\sqrt{r}}\|_2^2+\|\f{v_\theta}{\sqrt{r}}\|_2^2+\|\sqrt{r}Dv\|_2^2+\|\p_t v\|_{\tilde H^1}^2+\|\sqrt{r}{\p_t\phi}\|_2^2+\|\sqrt{r}\p_t^2\phi\|_2^2+\|\sqrt{r}\p_t^2v\|_2^2\nn\\
&&\quad +\|r^{1/2}\p_r({\bar\rho}^{\gamma-2}\phi)\|_2^2+\|r^{1/2} \p_r^2v_r\|_2^2+\|r^{1/2}\p_r({\bar\rho}^{\gamma-2}\p_z\phi)\|_2^2+\|r^{1/2}\p_r^2\p_z v_r\|_2^2+\|\sqrt{ r}\p_r^2({\bar\rho}^{\gamma-2}\phi)\|_2^2\biggl)d\tau\nn\\
&&\le C(\|\phi_0\|_{\tilde H^2}^2+\|v_0\|_{\tilde H^3}^2)+ C\int_0^t \biggl(A_1+A_2+A_3+\|  g\|_{\tilde H^1}^2+\| f\|_{\tilde H^1}^2
+|(r\p_r({\bar\rho}^{\gamma-2}\p_z\phi ), {\bar\rho}^{\gamma-2}\p_{rz}^2f)|\nn\\
&&\quad+|(\p_r^2({\bar\rho}^{\gamma-2}\phi),r{\bar\rho}^{\gamma-2}\p_r^2f)
| +|( \p_z^{2}g,{\bar\rho} r \p_z^{2}v)|+|( \p_z^{2}f,{\bar\rho}^{\gamma-2} r \p_z^{2}\phi)|\biggl)d\tau.
 \label{cc3.1}
\eea
Taking $k=0,j=1$ and $k=0,j=2$ in Lemma \ref{pcl3} respectively, and then adding them yields
\bea
&&   \|\sqrt{r}\p_z v\|_2^2+\|\sqrt{r}\p_z^2 v\|_2^2+\|\sqrt{r}\p_z \phi\|_2^2+\|\sqrt{r}\p_z^2 \phi\|_2^2
+\int_0^t \big(\|\sqrt{r}\p_z Dv\|_2^2+\|\sqrt{r}\p_z^2 Dv\|_2^2\big)d\tau\nn\\
&&\le C(\|\phi_0\|_{\tilde H^2}^2+\|v_0\|_{\tilde H^3}^2) +C \sum\limits_{j=1}^2\int_0^t \biggl(|(\p_z^jg, {\bar\rho} r \p^j_z v)|+|(\p_z^j f,{\bar\rho}^{\gamma-2} r\p_z^j \phi)|
 \biggl)d\tau.\label{cc3.2}
\eea
Set
\bea
M_1=A_1+A_2+A_3+\|  g\|_{\tilde H^1}^2+\| f\|_{\tilde H^1}^2+|(r D^2({\bar\rho}^{\gamma-2} \phi ), D^2({\bar\rho}^{\gamma-2}f))|
+|(D^2 g, {\bar\rho} r D^2 v)|. \label{cc3.3}
\eea
It follows from \eqref{cc3.1} and \eqref{cc3.2}   that
\bea
&&
\ds\|v\|_{\tilde H^1}^2+\|\p_t v\|_{\tilde H^1}^2+\|\sqrt{r}\phi\|_2^2+\|\sqrt{r}\p_t \phi\|_2^2+\|\sqrt{r}\p_z v\|_2^2+\|\sqrt{r}\p_z^2 v\|_2^2+\|\sqrt{r}\p_z \phi\|_2^2\nn\\
&& \quad \ds+\|\sqrt{r}\p_z^2 \phi\|_2^2
+
\|r^{1/2}\p_r({\bar\rho}^{\gamma-2}\phi)\|_2^2+\|r^{1/2}\p_r({\bar\rho}^{\gamma-2}\p_z\phi)\|_2^2+\|r^{1/2}\p_r^2({\bar\rho}^{\gamma-2}\phi)\|_2^2
\nn\\
&&\quad \ds+\int_0^t \biggl(\|\f{v_r}{\sqrt{r}}\|_2^2+\|\f{v_\theta}{\sqrt{r}}\|_2^2+\|\sqrt{r}Dv\|_2^2+\|\p_t v\|_{\tilde H^1}^2+\|\sqrt{r}{\p_t\phi}\|_2^2+\|\sqrt{r}\p_t^2\phi\|_2^2+\|\sqrt{r}\p_t^2v\|_2^2\nn\\
&&\quad \ds+\|r^{1/2}\p_r({\bar\rho}^{\gamma-2}\phi)\|_2^2+\|r^{1/2} \p_r^2v_r\|_2^2+\|r^{1/2}\p_r({\bar\rho}^{\gamma-2}\p_z\phi)\|_2^2+\|r^{1/2}\p_r^2\p_z v_r\|_2^2\nn\\
&&\quad \ds+\|\sqrt{ r}\p_r^2({\bar\rho}^{\gamma-2}\phi)\|_2^2+\|\sqrt{r}\p_z Dv\|_2^2+\|\sqrt{r}\p_z^2 Dv\|_2^2\biggl)d\tau \nn\\
&&\ds \le C\biggl(\|\phi_0\|_{\tilde H^2}^2+\|v_0\|_{\tilde H^3}^2+ \int_0^t M_1d\tau \biggl).
\label{cc3.4}
\eea
In addition, by (1.13), we see that
$$
D\p_t \phi =D(-\f{1}{r}\p_r(r {\bar\rho} v_r) -\p_z ({\bar\rho} v_z) +f).
$$
Rewriting \eqref{cc1.7} and \eqref{cc1.8} as follows
\bea
&&\nu_1 \p_r^2 v_\theta ={\bar\rho} \p_t v_\theta -\nu_1\p_r(\f{v_\theta}{r})-{\bar\rho}g_2, \label{cc3.5}\\
&&\nu_1 \p_r^2 v_z ={\bar\rho}\p_t v_z+\gamma {\bar\rho}^{\gamma-1}\p_z\phi -\nu_1\p_z^2v_z -\nu_1\f{1}{r}\p_r v_z -\nu_2 \p_z(\f{1}{r}\p_r(rv_r)+\p_z v_z)-{\bar\rho} g_3.\label{cc3.6}
\eea
Then by \eqref{cc3.4}-\eqref{cc3.6} and \eqref{cc2.34} with $j=0$ , we have
\bea
&& \|v\|_{\tilde H^2}^2+\|\p_t v\|_{\tilde H^1}^2+\|\rho^{\gamma-2}\phi\|_{\tilde H^2}^2+\|\p_t\phi\|_{\tilde H^1}^2 \nn\\
&&\le C\biggl(\|\phi_0\|_{\tilde H^2}^2+\|v_0\|_{\tilde H^3}^2+ \int_0^t M_1d\tau \biggl)+C\sup\limits_{0\le \tau \le t}(\|f\|_{\tilde H^1}^2+\|\sqrt{r}g\|_2^2). \label{cc3.7}
\eea
Furthermore, by   \eqref{cc3.5}-\eqref{cc3.6}, \eqref{cc2.34} with $j=0$,  \eqref{cc3.4} and \eqref{cc3.7},  we also
have
\bea
&& \|v\|_{\tilde H^3}^2+\|\p_t v\|_{\tilde H^1}^2+\|\rho^{\gamma-2}\phi\|_{\tilde H^2}^2+\|\p_t\phi\|_{\tilde H^1}^2 \nn\\
&&\quad +\int_0^t \biggl(\|\f{v_r}{\sqrt{r}}\|_2^2+\|\f{v_\theta}{\sqrt{r}}\|_2^2+\|\sqrt{r}Dv\|_{\tilde H^3}^2+\|\p_t v\|_{\tilde H^2}^2+\|\sqrt{r}{\p_t\phi}\|_{\tilde H^1}^2+\|\sqrt{r}\p_t^2\phi\|_2^2\nn\\
&&\quad +\|\sqrt{r}\p_t^2v\|_2^2+\|D({\bar\rho}^{\gamma-2}\phi)\|_{\tilde H^1}^2\biggl)d\tau \nn\\
&&\le C\biggl(\|\phi_0\|_{\tilde H^2}^2+\|v_0\|_{\tilde H^3}^2+ \int_0^t M_1d\tau \biggl)+C\sup\limits_{0\le \tau \le t}(\|f\|_{\tilde H^1}^2+\|g\|_{\tilde H^1}^2).\label{cc3.8}
\eea
Define the energy functional
\bea
N(t)&=&  \|v\|_{\tilde H^3}^2+\|\p_t v\|_{\tilde H^1}^2+\|\rho^{\gamma-2}\phi\|_{\tilde H^2}^2+\|\p_t\phi\|_{\tilde H^1}^2
+\int_0^t \biggl(\|\f{v_r}{\sqrt{r}}\|_2^2+\|\f{v_\theta}{\sqrt{r}}\|_2^2+\|\sqrt{r}Dv\|_{\tilde H^3}^2 \nn\\
&&+\|\p_t v\|_{\tilde H^2}^2+\|\sqrt{r}{\p_t\phi}\|_{\tilde H^1}^2+\|\sqrt{r}\p_t^2\phi\|_2^2+\|\sqrt{r}\p_t^2v\|_2^2+\|D({\bar\rho}^{\gamma-2}\phi)\|_{\tilde H^1}^2\biggl)d\tau.\nn
\eea
Next, we show that $N(t)$ is uniformly bounded for any $t\ge 0$. In the proof procedure,
we will employ the following Gagoliado-Nirenberg's inequality repeatedly:
\begin{lemma}\label{pcl3.1}
(i) Assume $2\le p\le +\infty$. Let $j$ and $k$ be integers satisfying
$$0\le j<k ,\  k>j+2(\f{1}{2}-\f{1}{p}).$$
Then there exists a constant $C>0$ such that
$$\|D^j w\|_{L^p(\Omega)}\le C \|w\|_{L^2(\Omega)}^{1-a}\|D^k w\|_{L^2(\Omega)}^a,$$
where $\ds a=\f{1}{k}(j+1-\f{2}{p})$.

(ii) Assume $w(x)=w(r,z)$ and $x\in\Bbb R^3$. Let $2\le p\le +\infty$ and let $j$ and $k$ be integers satisfying
$$0\le j<k , k>j+3(\f{1}{2}-\f{1}{p}).$$
Then there exists a constant $C>0$ such that
$$\|D^j w\|_{L_r^p(\Omega)}\le C \|w\|_{L_r^2(\Omega)}^{1-a}\|D^k w\|_{L_r^2(\Omega)}^a,$$
where $\ds a=\f{1}{k}(j+\f{3}{2}-\f{3}{p})$.
\end{lemma}
\begin{proposition}\label{pcp1} Assume  $N(t)\le 1$, we then have
\be
N(t)\le  C(\|v_0\|_{\tilde H^3}^2+\|\phi_0\|_{\tilde H^2}^2)+ CN(t)^{3/2}.\nn
\ee
\end{proposition}

\pf At first, we deal with $A_i$ $(i=1,2,3)$ in the expression of $M_1$ in (3.8).
In order to estimate the term $\ds\int_0^t A_1d\tau$, by the expression of $A_1$, we
require to treat $\ds\int_0^t (g_1,{\bar\rho} r v_r)d\tau$ and $\ds\int_0^t (f, {\bar\rho}^{\gamma-2} r\phi)d\tau$. Note that
\bea
&&\int_0^t (g_1,{\bar\rho} r v_r)d\tau\nn\\
&&= \int_0^t \int_\Omega \biggl(\f{v_\theta^2}{r}-v_r\p_r v_r -v_z\p_zv_r-\p_r Q({\bar\rho} ,\phi) -\f{\nu_1 \phi}{(\phi+{\bar\rho}){\bar\rho}}(\p_r(\f{1}{r}\p_r(rv_r))+\p_z^2 v_\th)\nn\\
&&\quad-\f{\nu_2 \phi}{(\phi+{\bar\rho}){\bar\rho}}\p_r(\f{1}{r}\p_r(rv_r)+\p_z v_r)\biggl) r{\bar\rho} v_rdrdzd\tau. \label{cc3.9}
\eea
Next, we treat each term in the right hand side of \eqref{cc3.9}. One has that
\bea
&& \int_0^t \int_\Omega v_\theta^2 {\bar\rho} v_r drdzd\tau\nn\\
&&\le C \int_0^t \int_\Omega |r^{1/3}v_\theta||\f{v_\theta}{r^{1/2}}||r^{1/6}v_r|drdzdt\nn\\
&&\le C \int_0^t \|r^{1/3}v_\theta\|_{L^2}\|\f{v_\theta}{r^{1/2}}\|_{L^3} \|v_r\|_{L_r^6}d\tau.\label{cc3.10}
\eea
By Lemma \ref{pcl3.1} (i) with $j=0, p=3$ and $k=1$, we have
$$
\|\f{v_\theta}{r^{1/2}}\|_{L^3} \le C \|\f{v_\theta}{r^{1/2}}\|_{L^2}^{2/3}\|D(\f{v_\theta}{r^{1/2}})\|_{L^2}^{1/3}.
$$
This, together with Lemma \ref{pcl3.1}(ii) with $j=0, p=6$ and $k=1$, we obtain that from \eqref{cc3.10}
\bea
&& \int_0^t \int_\Omega\f{v_\theta^2}{r} r{\bar\rho} v_r drdzd\tau\nn\\
&& \le C\sup\limits_{0\le \tau \le t} \|r^{1/2} v_\theta\|_2 \int_0^t\|\f{v_\theta}{r^{1/2}}\|_{L^2}^{2/3}\|D(\f{v_\theta}{r^{1/2}})\|_{L^2}^{1/3}\|Dv_r\|_{L_r^2} d\tau
\nn\\
&&\le CN(t)^{3/2}.\label{cc3.11}
\eea
And we also have that
\bea
&& \int_0^t \int_\Omega v_r\p_r v_r r{\bar\rho} v_r drdzd\tau\nn\\
&&\le C\int_0^t \|r^{1/3}\p_r v_r\|_2\|r^{1/2}v_r\|_2\|v_r\|_{L_r^6}d\tau\nn\\
&& \le C\sup\limits_{0\le \tau \le t} \|r^{1/2} v_r\|_2 \int_0^t\|r^{1/2}\p_r v_r\|_2\|Dv_r\|_{L_r^2} d\tau \nn\\
&&\le CN(t)^{3/2}\label{cc3.12}
\eea
and
\bea
&& \int_0^t \int_\Omega \p_r Q({\bar\rho} ,\phi) r{\bar\rho} v_r drdzd\tau\nn\\
&&\le C \int_0^t \int_\Omega ( |\phi \p_r ({\bar\rho}^{\gamma-2}\phi) rv_r| +|\phi^2 \p_r {\bar\rho} r v_r | )drdzd\tau \nn\\
&&\le C \int_0^t \|r^{1/3}\phi\|_3\|r^{1/2}\p_r ({\bar\rho}^{\gamma-2}\phi)\|_2 \|r^{1/6}v_r\|_{6}d\tau + \int_0^t\|\phi\|_{3}\|r^{1/6}{\bar\rho}^{\gamma-2}\phi\|_6 \|\f{v_r}{\sqrt{r}}\|_2 d\tau \nn\\
&&\le CN(t)^{3/2},\label{cc3.13}
\eea
here we point out that we have used the crucial fact of  $\p_r {\bar\rho} =O( r^{-3})$ in (3.13), and $\|r^{1/6} v_r\|_6=\|v_r\|_{L_r^6}\le C \|D v_r\|_{L^2_r} = C\|r^{1/2}Dv_r\|_2$ by Lemme \ref{pcl3.1}(ii) with $j=0, p=6$ and $k=1$.

On the other hand, it follows from direct computation that
\bea
&& \int_0^t \int_\Omega \f{\phi}{(\phi+{\bar\rho}){\bar\rho}}\biggl(\nu_1(\p_r(\f{1}{r}\p_r(rv_r))+\p_z^2v_r)
+\nu_2 \p_r(\f{1}{r}\p_r(rv_r)+\p_zv_z) \biggl)r{\bar\rho} v_r drdzd\tau\nn\\
&& \le  C \int_0^t \int_\Omega \biggl( |\phi \f{v_r}{\sqrt{r}}\f{v_r}{\sqrt{r}}|+|\phi D v_r v_r |+|\sqrt{r}\phi\sqrt{r}D^2 v_r v_r | \biggl)drd\tau\nn\\
&&\le CN(t)^{3/2}.\label{cc3.14}
\eea
Similarly to \eqref{cc3.13}, we arrive at
\bea
\int_0^t (f,{\bar\rho}^{\gamma-2} r\phi)d\tau &=& -\int_0^t \int_\Omega (\f{1}{r}\p_r (r \phi v_r)+\p_z(\phi v_z)){\bar\rho}^{\gamma-2} r\phi drdzd\tau\nn\\
&&=\int_0^t \int_\Omega (r \phi v_r\p_r ({\bar\rho}^{\gamma-2} \phi)+r \phi v_r\p_z ({\bar\rho}^{\gamma-2} \phi)) drdzd\tau\nn\\
&&\le C\int_0^t \int_\Omega( |r \phi^2 v_r\p_r{\bar\rho} |+|r \phi  v_rD({\bar\rho}^{\gamma-2} \phi)|) drdzd\tau\nn\\
&&\le CN(t)^{3/2}.\label{cc3.15}
\eea
Combining \eqref{cc3.9} and \eqref{cc3.11}-\eqref{cc3.15}, we eventually obtain
\be
\int_0^t A_1d\tau\le CN(t)^{3/2}.\label{cc3.16}
\ee
Analogously, $A_i$ ($i=2,3$) and the terms  such as $\|g\|_{\tilde H^1}^2$, $\|f\|_{\tilde H^1}^2$ and $(D^2g, {\bar\rho} r D^2v)$ can
be treated like $A_1$.
For examples, we treat the terms  $(\p_t g_3, {\bar\rho} r \p_t^2 v_z)$ and $(rD^2(\rho^{\gamma-2}\phi),D^2({\bar\rho}^{\gamma-2}f))$
in $A_3$. Note that
\bea
&&\int_0^t (\p_t g_3,{\bar\rho} r \p_{t}^2v_z)d\tau\nn\\
&&=\int_0^t \int_\Omega \p_t\biggl(-v_r\p_r v_z-v_z\p_z v_z  -\p_z Q({\bar\rho},\phi) -\f{\nu_1 \phi}{(\phi+{\bar\rho}){\bar\rho}}\biggl(\p_r^2 v_z+\p_z^2 v_z +\f{1}{r} \p_r v_r \biggl)\nn\\
&&\quad \quad -\f{\nu_2 \phi}{(\phi+{\bar\rho}){\bar\rho}}\p_z(\f{1}{r}\p_r (r  v_r)+\p_z v_z) \biggl){\bar\rho} r \p_{t}^2v_z drdzd\tau. \nn
\eea
Using Holder's inequality and expression of $N(t)$, we see that
\bea
&&\int_0^t \int_\Omega \p_t\biggl(-v_r\p_r v_z-v_z\p_z v_z \biggl){\bar\rho} r \p_{t}^2v_z drdzd\tau \nn\\
&& \le C \int_0^t(\|r^{1/2}\p_tv\|_2\| \p_rv_z\|_{\infty}+\|r^{1/2}\p_{rt}^2 v_z\|_2\| v\|_{\infty})\|r^{1/2}\p_{t}^2v_z\|_2d\tau\nn\\
&& \le C N(t)^{3/2}. \label{cc3.17}
\eea
By Lemma \ref{pcl3.1} (i) with $j=0, p=4$, we can obtain
\bea
&&\int_0^t \int_\Omega \p_{tz}^2 Q({\bar\rho},\phi){\bar\rho} r \p_{t}^2v_z drdzd\tau \nn\\
&& \le C \int_0^t(\|r^{1/2}\p_{tz}^2\phi\|_2\| \phi\|_{\infty}+\|r^{1/2}\p_t\phi\|_4\| \p_z \phi\|_{4})\|r^{1/2}\p_{t}^2v_z\|_2d\tau\nn\\
&&\le C (\| {\bar\rho}^{\gamma-2}\phi\|_{\tilde H^1}+\| {\bar\rho}^{\gamma-2}\p_t \phi\|_{\tilde H^1})
\int_0^t(\|r^{1/2}\p_t\phi\|_{\tilde H^1}+\|r^{1/2} \p_z \phi\|_{2})\|r^{1/2}\p_{t}^2v_z\|_2d\tau.\nn\\
&& \le C N(t)^{3/2}. \label{cc3.18}
\eea
By Lemma \ref{pcl3.1} (i) with $j=1, p=4$ and $j=0, p=4$, one has
\bea
&&\int_0^t \int_\Omega \p_t\biggl\{-\f{ \phi}{(\phi+{\bar\rho}){\bar\rho}}\biggl(\nu_1 (\p_r^2 v_z+\p_z^2 v_z +\f{1}{r} \p_r v_r )
 - \nu_2  \p_z(\f{1}{r}\p_r (r  v_r)+\p_z v_z)\biggl) \biggl\}{\bar\rho} r \p_{t}^2v_zdrdzd\tau \nn\\
 &&  \le C \int_0^t\biggl(\|\p_t \phi\|_4 \|Dv\|_{\tilde H^2}+\| \phi\|_{\infty} \|\p_tv\|_{\tilde H^2}\biggl)\|r^{1/2}\p_{t}^2v_z\|_2 d\tau \nn\\
 && \le C N(t)^{3/2}. \label{cc3.19}
\eea
Then collecting  \eqref{cc3.17}-\eqref{cc3.19} yields
\be
\int_0^t (\p_t g_3, {\bar\rho} r \p_{t}^2v_z)d\tau\le C N(t)^{3/2}. \label{cc3.20}
\ee
In the end, we deal with the term $\ds \int_0^t (rD^2(\rho^{\gamma-2}\phi),D^2({\bar\rho}^{\gamma-2}f))d\tau$. We see that
\bea
&& \int_0^t (rD^2(\rho^{\gamma-2}\phi),D^2({\bar\rho}^{\gamma-2}f))d\tau \nn\\
&&=\int_0^t \int_\Omega rD^2(\rho^{\gamma-2}\phi)D^2({\bar\rho}^{\gamma-2}(-\f{1}{r}\p_r (r\phi v_r)-\p_z(\phi v_z)))drdzd\tau.\label{cc3.21}
\eea
In addition, we observe that
\bea
&&D^2\biggl({\bar\rho}^{\gamma-2}(-\f{1}{r}\p_r (r\phi v_r)-\p_z(\phi v_z))\biggl)\nn\\
&&=-v_r \p_r(D^2(\rho^{\gamma-2}\phi))-v_z \p_z(D^2(\rho^{\gamma-2}\phi))  +\mathcal{R} \nn
\eea
with
$$
\mathcal{R} =D^2\biggl({\bar\rho}^{\gamma-2}(-\f{1}{r}\p_r (r\phi v_r)-\p_z(\phi v_z))\biggl)+v_r \p_r(D^2(\rho^{\gamma-2}\phi))+v_z \p_z(D^2(\rho^{\gamma-2}\phi)).
$$
It follows from direct computation that
\bea
&&\int_0^t \biggl|\int_\Omega rD^2(\rho^{\gamma-2}\phi) \mathcal{R}drdz\biggl|d\tau\nn\\
&&\le C (\|r^{1/2}D^2(\rho^{\gamma-2}\phi)\|_2  +\|v\|_{\tilde H^3})\int_0^t(\|r^{1/2}D^2(\rho^{\gamma-2}\phi)\|_2
+\|Dv\|_{\tilde H^3})\|Dv\|_{\tilde H^3}d\tau\nn\\
&&\le C N(t)^{3/2}. \label{cc3.22}
\eea
For the remainder term  in \eqref{cc3.21}, one sees that
\bea
&&\biggl|\int_0^t \int_\Omega rD^2(\rho^{\gamma-2}\phi)\biggl(-v_r \p_r(D^2(\rho^{\gamma-2}\phi))-v_z \p_z(D^2(\rho^{\gamma-2}\phi))\biggl)drdzd\tau\biggl|\nn\\
&&=\f{1}{2} \biggl|\int_0^t \int_\Omega (D^2(\rho^{\gamma-2}\phi))^2
(\p_r(r v_r)+r\p_z v_z)drdzd\tau\biggl|\nn\\
&& \le C N(t)^{3/2}. \quad \quad (\text{Analogously as in \eqref{cc3.22} }) \label{cc3.23}
\eea
Thus
\bea
&&\int_0^t (rD^2(\rho^{\gamma-2}\phi),D^2({\bar\rho}^{\gamma-2}f))d\tau \le C N(t)^{3/2}.\label{cc3.24}
\eea
Consequently, by \eqref{cc3.8}, \eqref{cc3.16}, \eqref{cc3.20}, \eqref{cc3.23}-\eqref{cc3.24}, and analogous treatments
for $A_1$, we obtain that
\bea
&&\sum \limits_{i=1}^3\bigl|\int_0^t A_id\tau\bigl|\le  CN(t)^{3/2}\label{pcr4.18}
\eea
and
\bea
&&\sup \limits_{0\le \tau\le t}(\|f\|_{\tilde H^1}^2+\|  g\|_{\tilde H^1}^2)
\le CN(t)^{3/2}.
\label{pcr4.19}
\eea
Substituting \eqref{pcr4.18}-\eqref{pcr4.19} into \eqref{cc3.8} yields the proof of Proposition 3.1.
\ef
\vskip 0.1 true cm

Based on Proposition 3.1, we now start to prove Theorem 1.1.

{\bf Proof of Theorem 1.1.}  By Proposition 3.1, when $\|v_0\|_{\tilde H^3}^2+\|\rho_0\|_{\tilde H^2}^2\le \ve^2 $ and
$\ve>0$ is small, then $N(t)\le C\ve^2$ holds uniformly for any $t\ge 0$. This, together
with the local existence of classical solution to (1.13)-(1.16) with (1.17)-(1.18) (one can see [20]) and continuity argument,
yields the global solution of  problem
(1.13)-(1.16) with (1.17)-(1.18). Thus the proof of Theorem 1.1 is completed.\ef


\begin{thebibliography}{llll}


\bibitem{1} Cho Yonggeun, Choe Hi Jun, Kim Hyunseok, {\it Unique solvability of the initial boundary value
problems for compressible viscous fluids},   J. Math. Pures Appl. 83, 243-275 (2004)

\bibitem{2} Choe Hi Jun, Kim Hyunseok, {\it Global existence of the radially symmetric solutions of
the Navier-Stokes equations for the isentropic compressible fluids}, Math. Methods Appl. Sci. 28,
no. 1, 1-28 (2005)


\bibitem{3} R. Courant, K. O. Friedrichs, {\it Supersonic flow and
shock waves}, Interscience Publishers Inc., New York, 1948.


\bibitem{4} Ding Shijin, Wen Huanyao, Yao Lei, Zhu Changjiang,
{\it Global spherically symmetric classical solution to compressible Navier-Stokes equations
with large initial data and vacuum}, SIAM J. Math. Anal. 44, no. 2, 1257-1278  (2012)

\bibitem{5} E. Feireisl, {\it Dynamics of Viscous Compressible Fluids}, Oxford Univ. Press, Oxford, 2004.


\bibitem{6} E. Feireisl, A. Novotn\'y, H. Petzeltov\'a, {\it On the existence of globally defined weak solutions
to the Navier-Stokes equations}, J. Math. Fluid Mech. 3, 358-392 (2001)



\bibitem{7}  D. Hoff, {\it Spherically symmetric solutions of the Navier-Stokes equations for compressible,
isothermal flow with large, discontinuous initial data}, Indiana Univ. Math. J., 41,
1225-1302  (1992)

\bibitem{ 8}  Huang Xiangdi,  Li Jing, Xin Zhouping, {\it Global well-posedness of classical solutions with
large oscillations and vacuum to the three-dimensional isentropic compressible
Navier-Stokes equations}, Comm. Pure Appl. Math. 65, 549-585 (2012)

\bibitem{9}  Y. Kagei, T. Kobayashi, {\it On large-time behavior of solutions to
the compressible Navier-Stokes equations in the half space in $\Bbb R^3$},
Arch. Ration. Mech. Anal. 165, no. 2, 89-159  (2002)


\bibitem{10}  Y. Kagei, T. Kobayashi, {\it Asymptotic behavior of solutions of the compressible
Navier-Stokes equations on the half space}, Arch. Ration. Mech. Anal. 177, no. 2, 231-330  (2005)



\bibitem{11} Y. Kagei, {\it Asymptotic behavior of solutions of the compressible Navier-Stokes
equation around the plane Couette flow}, J. Math. Fluid Mech. 13, no. 1, 1-31  (2011)

\bibitem{12}  Y. Kagei, {\it Global existence of solutions to the compressible Navier-Stokes equation
around parallel flows}, J. Differential Equations 251, no. 11, 3248-3295  (2011)


\bibitem{13} P. L. Lions, {\it Mathematical Topics in Fluid Dynamics}, Vol. 2, Compressible Models, Oxford
Science Publication, Oxford, 1998.



\bibitem{14} Jiang Song, {\it Global spherically symmetric solutions to the equations of a viscous polytropic ideal
gas in an exterior domain}, Comm. Math. Phys., 178, 339-374 (1996)



\bibitem{15}  Jiang Song, Zhang Ping, {\it Axisymmetric solutions of the 3D Navier-Stokes equations for compressible
isentropic fluids}, J. Math. Pures Appl. (9) 82, no. 8, 949-973  (2003)



\bibitem{16}  A. Matsumura, T. Nishida, {\it The initial value problem for the equations of motion of
viscous and heat-conductive gases}, J. Math. Kyoto Univ., 20, 67-104 (1980)

\bibitem{17} A. Matsumura, T. Nishida, {\it Initial-boundary value problems for the equations of motion of
compressible viscous and heat-conductive fluids}, Comm. Math. Phys., 89, 445-464 (1983)

\bibitem{18} R. Salvi, I. Straskraba, {\it Global existence for viscous compressible
uids and their behavior as $t\to\infty$},
J. Fac. Sci. Univ. Tokyo Sect. IA Math. 40, no. 1, 17-51  (1993)


\bibitem{19} Sun Wenjun, Jiang Song, Guo Zhenhua, {\it  Helically symmetric solutions to the
3-D Navier-Stokes equations for compressible isentropic fluids}, J. Differential Equations 222,
no. 2, 263-296  (2006)


\bibitem{20} A. Tani, {\it On the first initial-boundary value problem of compressible viscous fluid motion}, Publ.
Res. Inst. Math. Sci. Kyoto Univ., 13, 193-253 (1977)

\bibitem{21} Wen Huanyao, Zhu Changjiang, {\it Global symmetric classical and strong solutions of
the full compressible Navier-Stokes equations with
vacuum and large initial data}, arXiv:1109.5328v2 (2012)

\bibitem{22} Yin Huicheng, Zhang Lin, {\it The global stability of 2-D viscous axisymmetric
circulatory  flows}, Preprint, 2015 (to appear in Discrete Contin. Dyn. Syst. -A).




\end{thebibliography}
\end{document}